\documentclass[letterpaper, 10 pt, conference]{ieeeconf}  

\IEEEoverridecommandlockouts                              
\usepackage{amsmath} 
\usepackage{amssymb}  
\usepackage{graphicx} 
\usepackage{url}
\usepackage{cite}     
\usepackage{tikz}
\usetikzlibrary{shapes.geometric}
\usepackage{amsmath}
\usepackage{xcolor}
\usepackage{ulem}
\usepackage{orcidlink}
\usetikzlibrary{shapes, arrows.meta, positioning}

\newtheorem{assumption}{Assumption}

\newcommand{\ex}{{\rm ex}}

\newcommand{\e}{{\rm e}}
\newcommand{\n}{{\rm n}}
\newcommand{\g}{{\rm g}}
\newcommand{\s}{{\rm sc}}
\newcommand{\sg}{{\rm cg}}
\newcommand{\inp}{{\rm in}}
\newcommand{\out}{{\rm out}}

\newtheorem{remark}{Remark}

\DeclareMathOperator{\diag}{diag}

\title{\LARGE \bf 
A Port-Hamiltonian Modeling Approach for Integrated Hydrogen Systems
}

\author{Abdullah Shahin$^1$\,\orcidlink{0009-0008-4552-6934}, Hannes Gernandt$^{1,2}$\,\orcidlink{0000-0001-7364-4606}, Anton Plietzsch$^1$\,\orcidlink{0000-0001-7480-2493}, and Johannes Schiffer$^{1,3}$\,\orcidlink{0000-0001-5639-4326}  
\thanks{This work was supported by the European Union’s Horizon Europe Framework Programme (HORIZON) under the GA n. 101120278 - DENSE.
HG and JS acknowledge funding from the German Federal Government, the Federal Ministry of Education and Research, and the State of Brandenburg within the framework of the joint project EIZ: Energy Innovation Center (project numbers 85056897 and 03SF0693A)
}
\thanks{$^1$Fraunhofer IEG, Fraunhofer Research Institution for Energy Infrastructures and Geotechnologies IEG, 03046 Cottbus, Germany,
{\tt\small \{abdullah.shahin, anton.plietzsch\}@ieg.fraunhofer.de,}
}%
\thanks{$^{2}$Institute of Mathematical Modelling, Analysis and Computational Mathematics IMACM, University of Wuppertal, 42119 Wuppertal, Germany, {\tt\small gernandt@uni-wuppertal.de}}
\thanks{$^3$Brandenburg University of Technology Cottbus-Senftenberg, 03046 Cottbus, Germany, {\tt\small schiffer@b-tu.de}}
}
\newcommand{\iconstorage}{\tikz[baseline=-0.25em]{\node[diamond,draw=black,fill=green!30,inner sep=2pt]{};}}
\newcommand{\iconjunction}{\tikz[baseline=-0.25em]{\node[diamond,draw=black,fill=cyan!50,inner sep=2pt]{};}}
\newcommand{\iconcompressor}{\tikz[baseline=-0.25em]{\node[diamond,draw=black,fill=yellow!50,inner sep=2pt]{};}}
\newcommand{\iconelectrolyzer}{\tikz[baseline=-0.25em]{\node[ellipse,draw=black,fill=blue!30,inner sep=2pt]{};}}
\newcommand{\iconfuelcell}{\tikz[baseline=-0.25em]{\node[ellipse,draw=black,fill=red!30,inner sep=2pt]{};}}

\begin{document}
\floatsep=.99mm
\textfloatsep=.99mm
\abovedisplayskip=.99mm
\belowdisplayskip=.99mm

\maketitle
\thispagestyle{empty}
\pagestyle{empty}
	\begin{abstract}
	
	Hydrogen's growing role in the transition towards climate-neutral energy systems necessitates structured modeling frameworks. 
	Existing gas network models, largely developed for natural gas, fail to capture hydrogen systems distinct properties, particularly the coupling of hydrogen pipes with electrolyzers, fuel cells, and electrically driven compressors. In this work, we present a unified systematic port-Hamiltonian (pH) framework for modeling hydrogen systems,  
    which inherently provides a passive input-output map of the overall interconnected system and, thus, a promising foundation for structured analysis, control and optimization of this type of newly emerging energy systems.
\end{abstract}
\section{Introduction}

The shift to climate-neutral energy systems requires integrating various energy carriers, with hydrogen playing a crucial role in decarbonization \cite{NEUMANN20231793}. 
Yet, for hydrogen to indeed become a key energy carrier, a dedicated generation and transportation infrastructure along with suitable operation and control structures must be established \cite{YUE2021111180}. 
Such development is, e.g., already underway in Germany, where a ``hydrogen core network" is currently being implemented with planned commissioning in 2032 \cite{hydrogencorenetwork2024}. 

For climate-neutrality, the hydrogen production needs to be ``green", i.e., based on renewable energy \cite{DAWOOD20203847}. 
Likewise, hydrogen is also foreseen as an important carrier for providing long-term energy storage to electric power systems \cite{YUE2021111180}. This results in a bidirectional energy conversion between electrical and hydrogen domains, so that these systems exhibit complex interdependencies that require structured, modular modeling approaches. Hence, hydrogen network components, which include production, storage, and distribution, must interact seamlessly with electrical grids to facilitate efficient and climate-neutral energy provision.

Clearly, many new control and operational challenges arise for these type of networks. A first step towards addressing these challenges in a structured manner is the development of suitable modeling procedures for systematically assembling mathematical models of hydrogen systems that reliably capture their fundamental dynamics. This is the main objective pursued in the present paper.

In this context, port-Hamiltonian (pH) systems offer a powerful framework to represent energy-based multi-domain systems, while preserving structural and systemic properties, such as interconnection patterns and passivity \cite{porthamiltonian_schaft_2014}. The pH framework has been used successfully to model electric power \cite{fiaz2013port,schiffer2014conditions,chargingstationpaper}, heat \cite{hauschild2020port,strehle2022port,krishna2021port} and gas networks \cite{domschke2023network}.

In the gas network domain, recently, in \cite{malan2023port} a pH representation of lumped gas pipe models and their interconnection has been developed. But therein the focus lies on pipeline dynamics for natural gas systems, omitting hydrogen storage, generation, and load integration, such as electrolyzers and fuel cells. Building on \cite{malan2023port}, in \cite{MalanGiesslerStrehle2024} compressors were introduced into the gas grid model. In \cite{MalanGiesslerStrehle2024} the control theoretic properties of a lumped compressor-gas network were analyzed through the establishment of equilibrium-independent passivity (EIP). Additionally, a pH PDE formulation of the gas network including compressors and based on the Euler equations was derived in \cite{Bendokat24}. Their model explicitly accounts for the enthalpy added to the gas in the compression process, providing a detailed representation of the involved thermodynamics.

Only a few recent studies have considered hydrogen-electric systems within the pH framework. The existing works \cite{KUMAR2024111814} and \cite{SBARBARO201838} primarily focus on the component-level modeling of PEM fuel cells and electrochemical processes, offering insight into the internal dynamics and energy conversion behavior of these individual units. However, these models are not readily extendable to network-scale applications, where the interaction between multiple components and energy carriers must be coherently represented. 

In light of this, the main contribution of our work is to adopt a system-level perspective and systematically develop a unified pH model that captures the dynamics of interconnected hydrogen infrastructures. To this end, we first derive individual pH representations for the components of the core hydrogen system, such as pipes, electrolyzers, fuel cells, and electrically driven compressors.
Next, by using algebraic graph theory and Kirchhoff's law to describe the network interconnections, we provide an overall pH model, which ensures structure-preserving interconnections among system components. 
In this way, the presented approach
reveals that hydrogen systems admit a compositional pH structure, which inherently provides a passive input-output map of the overall interconnected system and, thus, a promising
foundation for structured analysis, control, and optimization of this type of newly emerging energy systems.

The remainder of this paper is structured as follows: In Section~\ref{sec:monotone_phs}, we recall the class of pH systems. The pH models of the hydrogen grid components are presented in Section~\ref{sec:pH_components}. In Section~\ref{sec:network} the interconnection and integration of hydrogen components in an interconnected system is described.

\textbf{Notation.} We employ the following short-hands.
For $\ell>1$ scalars $v_i,$ $i=1,\ldots,\ell,$ $v=(v_{i})_{i=1}^\ell\in\mathbb R^\ell$ denotes a column vector, whose entries are $v_i.$ For $v\in\mathbb{R}^\ell$ and $w\in\mathbb{R}^{\ell'}$, where $\ell'>1$, we denote $(v,w)\in\mathbb{R}^{\ell+\ell'}$ as a column vector containing the entries of $v$ and $w$. 
Given $\ell>1$ matrices $M_i\in\mathbb{R}^{n_i\times m_i}$, $i=1,\ldots,\ell,$ the notation $\diag(M_i)_{i=1}^\ell \in\mathbb{R}^{(\sum_{i=1}^\ell n_i)\times (\sum_{i=1}^\ell m_i)}$ denotes the block matrix, whose block elements are $M_i$. The identity matrix of size $k\times k$ is denoted by $I_k$. Furthermore, the zero matrix of size $k\times l$ is denoted by $0_{k\times l}$ and in the case $l=1$ we simply write $0_k\in\mathbb{R}^k$.


\section{Port-Hamiltonian systems}
\label{sec:monotone_phs}
In the following, we recall from \cite{vanderSchaft2017} the form of pH systems employed throughout the paper 
\begin{align}
\label{eq:usual_phs}
\begin{split}
\dot x&=(J(x)-R(x))\nabla H(x)+Bu+d,\\
y&=B^\top\nabla H(x)+Du,
\end{split}
\end{align}
where $H:\mathbb{R}^n\rightarrow[0,\infty)$ is the Hamiltonian, which is assumed to be continuously differentiable with gradient $\nabla H(x)$, $B\in\mathbb{R}^{n\times m}$ is the input matrix, $D\in\mathbb{R}^{m\times m}$ is the feedthrough matrix and $d\in\mathbb{R}^n$ is a constant disturbance. Furthermore, the interconnection and damping matrices $J:\mathbb{R}^n \to \mathbb{R}^{n\times n}$ and $R:\mathbb{R}^n \to \mathbb{R}^{n\times n}$, respectively, satisfy the following structural conditions for all $x\in\mathbb{R}^n,$
\begin{align*}
J(x)=-J(x)^\top\quad R(x)=R(x)^\top\geq 0,\quad 
D+D^\top \geq 0.
\end{align*}
This directly implies that in the unperturbed case, i.e. $d\equiv 0,$ the system \eqref{eq:usual_phs} is passive \cite{vanderSchaft2017}, i.e.\ 
the time derivative of the Hamiltonian $H$ along the solutions of the system \eqref{eq:usual_phs} satisfies
\begin{equation}
\begin{split}
\dot H&=-\nabla H^\top R(x)\nabla H+\nabla H^\top Bu\\
&=-\nabla H^\top R(x)\nabla H+y^\top u -\tfrac12u^\top(D+D^\top)u\leq y^\top u.    
\end{split}
\notag
\end{equation}
 
\section{Port-Hamiltonian modeling of hydrogen systems}
\label{sec:pH_components}
In this section, we introduce pH models for pipes, compressors and storage units and also develop pH representations for electrolyzers and fuel cells.

\subsection{Network topology}
\label{sec:topol}
The hydrogen network topology is described by an undirected, connected graph $\mathcal{G}=(\mathcal{N},\mathcal{E})$ with set of nodes $\mathcal{N}=\{\n_1,\ldots,\n_N\}$ and set of edges $\mathcal{E}=\{\e_1,\ldots,\e_M\}$. To each node $\n_i\in\mathcal{N}$, we associate a node pressure $\mathbf{p}_i\in\mathbb{R}$. Moreover, to each edge $\e_\ell=\{\n_i,\n_k\}\in\mathcal{E}$, 
we associate a volumetric flow rate $\mathbf{q}_{\ell}\in\mathbb R$ and assign an arbitrary direction to each edge. 
Then, if $\mathbf{q}_\ell$ is directed from $\n_i$ to $\n_k$, we call $\n_i$ the source of $\e_\ell$ and $\n_k$ the sink of $\e_\ell$ and we write $\e_\ell=\overrightarrow{\n_i\n_k}$. 
It is convenient to introduce the node-edge incidence matrix $B_{\mathcal{G}}=(b_{ij}) \in \mathbb{R}^{N \times M}$, which is defined by
\begin{align}
b_{ij} =
\begin{cases}
1 & \text{if node  $\n_i$ is the sink of edge $\e_j$}, \\
-1 & \text{if node $\n_i$ is the source of edge $\e_j$,} \\
0 & \text{otherwise.}
\end{cases} 
\end{align}

We consider a hydrogen system that contains $S\geq 1$ storage units, $F\geq 1$ junctions, and $C\geq 1$ compressors. The set of nodes $\mathcal N$ is thus subdivided into the set $\mathcal{N}_S=\{\n_1,\ldots,\n_S\}$ that represents storage units, the set $\mathcal{N}_F=\{\n_{S+1},\ldots,\n_{S+F}\}$ that represents junctions within the hydrogen grid and the set of compressors $\mathcal{N}_C=\{\n_{S+F+1},\ldots,\n_{S+F+C}\},$ with $\n_{S+F+C}=\n_N$ and $\mathcal{N}=\mathcal{N}_S\cup\mathcal{N}_F\cup\mathcal{N}_C.$ 

The hydrogen system is coupled to the electric domain via $E\geq0$ electrolyzers and $L\geq0$ fuel cells. The interconnection topology between the electric power system and the hydrogen system is modeled by a second undirected graph $\mathcal G_P=({\mathcal{N}_S}\cup\mathcal{N}_P,\mathcal E_P).$
As commonly done in practice, we assume that each power conversion device is always connected to a dedicated storage unit on the hydrogen side, i.e., to a node $\n_i\in\mathcal N_S.$ Hence, $E+L\leq S$. Consequently, the corresponding buses in the electrical power system are represented by the set $\mathcal{N_P}=\{\n_{P,1},\ldots,\n_{P,E+L}\}$.
To each $\n_{P,i}\in\mathcal N_P$, 
we associate a voltage $\mathbf{v}_i\in\mathbb R$. 
Furthermore, each pair of nodes $\{\n_i,\n_{P,i}\}$ is interconnected by an edge $\e_{P,k}=\{\n_i,\n_{P,i}\}\in\mathcal E_P$ with $\mathcal E_P=\{\e_{P,1},\ldots,\e_{P,E+L}\},$ where the edge $\e_{P,k}\in\mathcal E_P$ represents the volumetric flow rate between a fuel cell or an electrolyzer and the hydrogen grid, respectively. 

In the sequel, we focus on developing the main components of a sector coupled hydrogen system. 
The vector of pressures of the hydrogen system is denoted by $\mathbf{p}=(\mathbf{p}_i)_{i=1}^N$ and that of the volumetric flow rates by $\mathbf{q}=(\mathbf{q}_i)_{i=1}^M.$ Moreover, we introduce the vector of exogenous volumetric flow rates $\mathbf{q}_\ex=(\mathbf{q}_{\ex,i})_{i=1}^{S+F}\in\mathbb{R}^{S+F}$, 
that is used to represent exogenous hydrogen injections or demands stemming, e.g., from onshore terminals. 

\begin{figure}[htbp]
\centering
\begin{minipage}{\columnwidth}  
\centering
\begin{tikzpicture}[scale=0.75, every node/.style={scale=0.75}, node distance=0.75cm]

\tikzstyle{junction}=[diamond, draw=black, fill=cyan!50, minimum size=5.5mm]
\tikzstyle{storage}=[diamond, draw=black, fill=green!30, minimum size=5.5mm]
\tikzstyle{comp}=[diamond, draw=black, fill=yellow!50, minimum size=5.5mm]
\tikzstyle{ely}=[ellipse, draw=black, fill=blue!30, minimum size=5.5mm]
\tikzstyle{fc}=[ellipse, draw=black, fill=red!30, minimum size=5.5mm]

\node[junction] (j1) {$\n_{4}$};
\node[storage, above=1cm of j1] (s1) {$\n_1$};

\node[ely, above=1cm of s1] (ely1) {$\n_{P,1}$};
\node[junction, right=1cm of j1] (j2) {$\n_{5}$};

\node[storage, above=1cm of j2] (s2) {$\n_2$};

\node[fc, above=1cm of s2] (fc1) {$\n_{P,2}$};

\node[junction, right=1cm of j2] (j3) {$\n_{6}$};
\node[comp, above=1cm of j3] (c1) {$\n_{7}$};
\node[storage, right=1cm of c1] (s3) {$\n_{3}$};

\draw[-{Stealth}] (s1)--(j1) node[midway, left] {$\e_1$};
\draw[-{Stealth}] (s1)--(j1) node[midway, right] {$\mathbf{q}_{1}$};
\draw[-{Stealth}] (ely1)--(s1) node[midway, left] {$\e_{P,1}$};
\draw[-{Stealth}] (ely1)--(s1) node[midway, right] {$\mathbf{q}_{\mathrm{sc},1}$};
\draw[-{Stealth}] (j2)--(s2)node[midway, left] { $\e_3$};
\draw[-{Stealth}] (j2)--(s2)node[midway, right] {  $\mathbf{q}_{3}$};
\draw[-{Stealth}] (s2)--(fc1) node[midway, left] { $\e_{P,2}$ };
\draw[-{Stealth}] (s2)--(fc1) node[midway, right] {$\mathbf{q}_{\mathrm{sc},2}$};

\draw[-{Stealth}] (j1)--(j2)node[midway, above] {$\e_2$};
\draw[-{Stealth}] (j1)--(j2)node[midway, below] {$\mathbf{q}_{2}$};
\draw[-{Stealth}] (j2)--(j3)node[midway, above] {$\e_4$};
\draw[-{Stealth}] (j2)--(j3)node[midway, below] {$\mathbf{q}_{4}$};

\draw[-{Stealth}] (j3)--(c1) node[midway, left] {$\e_5$};
\draw[-{Stealth}] (j3)--(c1) node[midway, right] {$\mathbf{q}_{5}$};
\draw[-{Stealth}] (c1)--(s3) node[midway, above] {$\e_6$};
\draw[-{Stealth}] (c1)--(s3) node[midway, below] {$\mathbf{q}_{6}$};
\node[above=1cm of s3] (ex3) {};
\draw[dashed,-{Stealth}] (s3) -- (ex3)
  node[midway, right] {$\mathbf{q}_{\mathrm{ex},3}$};

\node[below=-0.1cm of j1] (ex41) {};
\node[left=1cm of ex41] (ex4) {};
\draw[dashed,-{Stealth}] (ex41) -- (ex4)
  node[midway, below] {$\mathbf{q}_{\mathrm{ex},4}$};

\node[below=-0.1cm of j2] (ex51) {};
\node[right=1cm of ex51] (ex5) {};
\draw[dashed,-{Stealth}] (ex51) -- (ex5)
  node[midway, below] {$\mathbf{q}_{\mathrm{ex},5}$};

\node[below=-0.1cm of j3] (ex61) {};
\node[right=1cm of ex61] (ex6) {};
\draw[dashed,-{Stealth}] (ex61) -- (ex6)
  node[midway, below] {$\mathbf{q}_{\mathrm{ex},6}$};

\end{tikzpicture}

\vspace{2mm}

\caption{%
Example hydrogen network topology.\, The symbols denote the following elements:
\protect \iconstorage\ storage unit $\mathcal{N}_{S}$,\,
\protect\iconjunction\ junction $\mathcal{N}_{F}$,\,
\protect\iconcompressor\ compressor $\mathcal{N}_{C}$,\,
\protect\iconelectrolyzer\ electrolyzer $\mathcal{N}_{P}$,\,
\protect\iconfuelcell\ fuel cell $\mathcal{N}_{P}$.\,
Arrows indicate the volumetric flow rates in pipelines $\mathbf{q}_{1},\mathbf{q}_{2},\mathbf{q}_{3},\mathbf{q}_{4}$, compressors $\mathbf{q}_{5}$ and $\mathbf{q}_{6}$, electrolyzer outlet $\mathbf{q}_{\mathrm{sc},1}$, fuel cells inlet $\mathbf{q}_{\mathrm{sc},2}$ and exogenous volumetric flow rates $\mathbf{q}_{\mathrm{ex},i}$. 
}

\label{fig:HydrogenNetwork}
\end{minipage}
\end{figure}
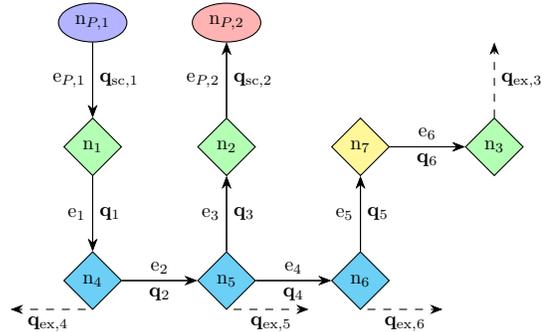

\subsection{Pipe model}
The hydrogen pipe associated with the $i$th edge $\e_i=\{\n_l,\n_r\}\in\mathcal E$, $\n_l\in\mathcal{N}_S\cup\mathcal{N}_F,$ $\n_r\in\mathcal{N}_S\cup\mathcal{N}_F,$ is considered to be a spatially one-dimensional object and the temporal change of pressure and volumetric flow rate of the hydrogen along the pipe can be described in terms of the isothermal compressible Euler equation, see e.g.\ \cite{domschke2023network}. 
A simplified lumped model which models the volumetric flow rate~$\mathbf{q}_{i}$ at standard conditions and with a constant compressibility factor was described in \cite[Proposition 2]{malan2023port} and is given by 
\begin{equation}
\label{eq:pipeflow}
   \rho \frac{\dot{\mathbf{q}}_{i}}{A_i}  = \frac{\mathbf{p}_l - \mathbf{p}_r}{L_i} 
    - \hat\lambda_i |\mathbf{q}_{i}| \mathbf{q}_{i}
    - \frac{g \sin(\theta_i)}{c^2} \textbf{p}_{M,i},
\end{equation}
with pressures $\mathbf{p}_l$ and $\mathbf{p}_r$ at the left and right ends of the pipeline segment, respectively,  the cross-sectorial area $A_i>0$, the diameter $D_i>0$, the pipe-length $L_i>0$, the speed of sound in hydrogen $c>0$, hydrogen density at standard conditions $\rho>0$, a friction coefficient 
\[
\hat\lambda_{i}=\frac{\lambda c^2 \rho^2 }{2 D_i A_i^2 \textbf{p}_{M,i}},
\]
with Darcy friction factor $\lambda>0$, pipe inclination angle $\theta_i\in[-\pi/2,\pi/2]$ and the gravitational acceleration $g>0$.

Moreover, $\mathbf{p}_{M,i}>0$ is the mean pressure across the $i$th pipeline segment given by the Weymouth mean pressure
\begin{equation}
\label{mean pressure}
    \mathbf{p}_{M,i} = \frac{2}{3} \frac{\mathbf{p}_l^3 - \mathbf{p}_r^3}{\mathbf{p}_l^2 - \mathbf{p}_r^2}
    = \frac{2}{3} \left( \mathbf{p}_l + \mathbf{p}_r - \frac{\mathbf{p}_l \mathbf{p}_r}{\mathbf{p}_l + \mathbf{p}_r} \right).
\end{equation}

To model the hydrogen flow rate within a pipe, we consider the following standard modeling assumptions, see also \cite{malan2023port}. 
\begin{assumption}
In the model \eqref{eq:pipeflow}, the inclination angle $\theta_i$ and the mean pressure $\textbf{p}_{M,i}$ given by \eqref{mean pressure} are constant.
 \label{ass:p}
\end{assumption}

With Assumption~\ref{ass:p} and by introducing the states, Hamiltonian and co-states 
\begin{align*}
x_{\e,i} = \rho\frac{L_i}{A_i} \mathbf{q}_i,\,
H_{\e,i}(x_{\e,i}) =  \tfrac{A_i}{2 L_i\rho}\|x_{\e,i}\|^2,\, 
\nabla H_{\e,i}(x_{\e,i}) = \mathbf{q}_i,
\end{align*} 
the pipe model~\eqref{eq:pipeflow} can be cast in the following pH form, see also~\cite[Theorem 3]{malan2023port},
\begin{align}
\dot{x}_{\e,i} &= (J_{\e,i}-R_{\e,i}(x_{\e,i}))\nabla H_{\e,i}(x_{\e,i})+ B_{\e,i}u_{\e,i} +d_{\e,i},\nonumber \\
y_{\e,i} &= B_{{\rm e},i}^\top \nabla H_{\e,i}(x_{\e,i})= \mathbf{q}_i,\quad u_{\e,i}=\mathbf{p}_l-\mathbf{p}_r, \label{eq:pipe_pH}\\
d_{\e,i} &=-\frac{gL_i \sin(\theta_i)}{c^2}\mathbf{p}_{M,i},\nonumber \\
B_{\e,i}&=1
\quad J_{\e,i}=0,\quad
 R_{\e,i}(x_{\e,i})= \hat\lambda_{i}L_i \frac{A_i}{L_i\rho}\big|x_{\e,i}\big|.  \nonumber 
\end{align}

As shown in \cite{malan2023port}, the lumped model \eqref{eq:pipe_pH} achieves an accuracy comparable to advanced discretized pipeline models~\cite{Pambour16}. Also, with Assumption~\ref{ass:p}, $R_{\e,i}(x_{\e,i})\geq 0$ holds for all $x_{\e,i}\in\mathbb{R}$. In comparison with \cite{malan2023port}, we do not incorporate the pressure dynamics of the boundary nodes of the pipe directly into the model \eqref{eq:pipeflow}. Instead, we incorporate these pressure dynamics in the junction and storage unit models, since this simplifies the modular interconnection of the individual network components.

\subsection{Hydrogen storage unit model}
In this section, we describe a general model of a hydrogen storage unit connected at node $\n_i\in\mathcal N_S$ with pressure $\mathbf{p}_{i}$ and which is interconnected via two edges $e_l\in\mathcal E$ and $e_k\in\mathcal E$.
Then, the dynamics is given by \cite{Klopcic24}  
\begin{equation}
	\frac{M_{H_2}V_{s,i}}{\rho RT_{s,i}}\dot{\mathbf{p}}_{i}=-\rho^{-1}r_{s,i}\mathbf{p}_{i}+\mathbf{q}_{\inp,i}-\mathbf{q}_{\out,i} +\mathbf{q}_{\ex,i},
			\label{eq:storage}
		\end{equation}
		where $V_{s,i}>0$ is the storage unit volume, $M_{H_2}>0$ is the molar mass of the stored hydrogen, $\mathbf{q}_{\inp,i}$ is the ingoing volumetric flow rate from the grid to the node $\n_i$, $\mathbf{q}_{\out,i}$ is the outgoing volumetric flow rate from the node $\n_i$ into the grid, 
        $\mathbf{q}_{\ex,i}$ is a volumetric flow rate that represents exogenous hydrogen injections or demands, $T_{{s,i}}>0$ is the temperature of the stored hydrogen, which is assumed to be constant, and $R=8.314 \ {\rm J}/({\rm mol}\cdot {\rm K})$ is the universal gas constant. Moreover, to account for potential storage unit losses, such as leakage, we introduce a dissipation constant $r_{s,i}>0$.

To write the dynamics \eqref{eq:storage} in pH form, we introduce 
\begin{equation*}
\begin{split}
 x_{\n,i}&=\frac{M_{H_2}V_{{s,i}}}{\rho RT_{{s,i}}}\mathbf{p}_{i},\quad H_{\n,i}(x_{\n,i})=\frac12\frac{\rho RT_{s,i}}{M_{H_2}V_{s,i}}x_{\n,i}^2,\\
R_{\n,i}&= r_{s,i}\geq 0, \quad 
\quad J_{\n,i}=0, \quad B_{\n,i}=\begin{bmatrix}1 & 1 \end{bmatrix},
\end{split} 
\end{equation*}
and define the output as well as the input
$$ 
y_{\n,i}=B_{\n,i}^\top\nabla H_{\n,i}(x_{\n,i})=\begin{bmatrix}
\mathbf{p}_{i}\\\mathbf{p}_{i}
\end{bmatrix},  \quad u_{\n,i}=\begin{bmatrix}
    \mathbf{q}_{\inp,i}-\mathbf{q}_{\out,i}\\ \mathbf{q}_{\ex,i}
\end{bmatrix}.
$$ 
With this, we obtain from \eqref{eq:storage} the following pH model
\begin{equation}
\label{eq:storage_pH}
\begin{split}
\dot x_{\n,i}&=-R_{\n,i}\nabla H_{\n,i}(x_{\n,i})+B_{\n,i}u_{\n,i},\\
y_{\n,i}&=B_{\n,i}^\top\nabla H_{\n,i}(x_{\n,i}).
\end{split}
\end{equation}

\subsection{Junction model}

Based on the considerations in \cite{malan2023port} we describe a  junction between different pipes at a node $n_i\in\mathcal{N}_F$ by the following equation for the node pressure 

\begin{align*}
C_i\dot{\mathbf{p}}_i=\mathbf{q}_{\inp,i}-\mathbf{q}_{\out,i} +\mathbf{q}_{\ex,i},\quad C_{i} = \sum_{l=1, \n_i\in\e_l}^M\frac{L_{l} A_{l}}{2 \rho c^2},
\end{align*}
where $C_i>0$ models a lumped storage capacity at the junction node $\n_i$ and 
we sum in the definition of $C_i$ over all edges that are incident with the node $\n_i$. 
Our junction model can be viewed as a lossless storage with a comparably small storage capacity $C_i$ and accordingly, the pH formulation of junctions is given similar to \eqref{eq:storage_pH} with  
\begin{equation}
\label{eq:junction_pH}
\begin{split}
 x_{\n,i}&=C_i\mathbf{p}_i,\quad H_{\n,i}(x_{\n,i})=\frac12C_i^{-1}x_{\n,i}^2,\\
R_{\n,i}&= 0, \quad 
\quad J_{\n,i}=0, \quad B_{\n,i}=\begin{bmatrix}1 & 1 \end{bmatrix}.
\end{split}
\end{equation}

\subsection{Compressor unit model}

We consider a centrifugal compressor at node $\n_i\in\mathcal N_C$, which is modeled according to \cite{Gravdahl97}. 
That is, we describe the pressure dynamics within the plenum $\mathbf{p}_i$ and using the input pressure $\mathbf{p}_{l}$ at the inlet grid node $\n_l\in{\mathcal{N}_S}\cup\mathcal{N}_F$ and the pressure $\mathbf{p}_{r}$ at the outlet  node $\n_r\in{\mathcal{N}_S}\cup\mathcal{N}_F$ that receives the volumetric flow rate from the plenum $\mathbf{p}_i$. This means that the corresponding volumetric flow rates $\mathbf{q}_{f}$ towards the plenum and from the plenum to the grid, $\mathbf{q}_{m}$, are represented via the volumetric flow rates over the edges $\e_f=\{\n_l,\n_i\}\in\mathcal{E}$ and $\e_m=\{\n_i,\n_r\}\in\mathcal{E}$. 
Thus, the dynamics are 
\begin{align}
\label{eq:comp_p_i}
\begin{split}
\frac{V_{p,i}}{\rho a_{01,i}^2}\dot{\mathbf{p}}_{i}&=-\rho^{-1}r_{pl,i}\mathbf{p}_{i}+\mathbf{q}_{f}-\mathbf{q}_{m}, \\
\frac{ \rho L_{c,f}}{A_{1,f}}\dot{\mathbf{q}}_{f}&=\mathbf{p}_{2,f}-\mathbf{p}_{i}, \\ 
\frac{ \rho L_{o,m}}{A_{o,m}} \dot{\mathbf{q}}_{m}&= \mathbf{p}_{i} - \mathbf{p}_{r}, 
\end{split}
\end{align}
where $a_{01,i}>0$ is the inlet stagnation sonic velocity, $V_{p,i}>0$ is the volume of the plenum, $L_{c,f}>0$ is the length of compressor and duct, $L_{o,m}>0$ the length of the compressor outlet, $A_{1,f}>0$, $A_{o,m}>0$ are the reference areas, $\mathbf{p}_{2,f}$ is the output pressure of the compressor and $r_{pl,i}>0$ is a coefficient accounting for pressure losses in the plenum. 

Adherent to the methodology proposed in \cite{Bendokat24}, the differential pressure across the compressor $\Delta \mathbf{p}_j$ is used as a control input, i.e., 
\begin{align}
\label{eq:comp2}
\mathbf{p}_{2,f}=\mathbf{p}_{2,f}-\mathbf{p}_{l}+\mathbf{p}_{l},\quad \Delta \mathbf{p}_{i-N+C}:=\mathbf{p}_{2,f}-\mathbf{p}_{l},
\end{align}
where the node index $i$ fulfills $i\in\{N-C+1,\ldots,N\}$ and therefore $i-N+C\in\{1,\ldots,C\}$.

In the following, we split the compressor model into three pH submodels: the pressure node dynamics in the plenum which can be modeled as a storage node \eqref{eq:storage_pH} with
\begin{align}
\nonumber
x_{\n,i}&=\frac{V_{p,i}}{\rho a_{01,i}^2}\mathbf{p}_{i},\quad  H_{\n,i}(x_{\n,i})=\frac12 \frac{\rho a_{01,i}^2}{V_{p,i}}x_{\n,i}^2,\quad d_{\n,i}=0,\\ R_{\n,i}&=\frac{r_{pl,i}}{\rho}, \quad J_{\n,i}=0,\quad  B_{\n,i}=1,\quad  
u_{\n,i}=\mathbf{q}_f-\mathbf{q}_m.\label{eq:comp_plenum_pH}
\end{align}
Secondly, the volumetric flow rate through the throttle can be modeled analogously to the pH pipe model \eqref{eq:pipe_pH}
with
\begin{align} 
\nonumber
x_{\e,m}&=\frac{ \rho L_{o,m}}{A_{o,m}} \mathbf{q}_{m},\,  H_{\e,m}(x_{t,m})=\frac12 \frac{ A_{o,m}}{\rho L_{o,m}}  x_{\e,m}^2,\, d_{\e,m}=0,\\ R_{\e,m}&=J_{\e,m}=0,\quad  B_{\e,m}=1,\quad
u_{\e,m}=\mathbf{p}_i-\mathbf{p}_r.
\label{eq:comp_throttle_pH}
\end{align}

Similarly, the volumetric flow rate dynamics through the compressor towards the plenum can be modeled analogously to the pH pipe model \eqref{eq:pipe_pH} with 
\begin{align}
\label{eq:comp_flow_pH}
x_{\e,f}&=\frac{ \rho L_{c,f}}{A_{1,f}} \mathbf{q}_{f},\,  H_{\e,f}(x_{\e,f})=\frac{ A_{1,f}}{2\rho L_{c,f}}  x_{\e,f}^2,\, d_{\e,f}=0, \\ R_{\e,f}&=J_{\e,f}=0,\quad  
B_{\e,f}=\begin{bmatrix}1 & 1\end{bmatrix}
,\quad  u_{\e,f}=\begin{bmatrix}
    \mathbf{p}_l-\mathbf{p}_i\\ \Delta \mathbf{p}_f 
\end{bmatrix}. \nonumber
\end{align}

\begin{remark}
Following \cite{Gravdahl97} one can also impose closed coupled valve control to replace the pressure difference $\Delta\mathbf{p}_{f}$ in \eqref{eq:comp2} by a desired dissipating term. Additionally, the model could be augmented by incorporating the compressor shaft's rotational speed $\omega$ as a state variable, with compressor torque as the control input.
\end{remark}

\subsection{Sector coupling components}

In this section, we present lumped pH electrolyzer and fuel cell models that facilitate the integration of the hydrogen system with the electrical domain that can be used for control from a macroscopic grid perspective of an integrated hydrogen network. Recall from Section~\ref{sec:topol} that the interconnection topology between the hydrogen system and the electric power system is modeled by the graph $\mathcal G_P=(\mathcal{N}_F\cup\mathcal{N}_P,\mathcal E_P)$.

First, we describe the electrolyzer model. Due to the increased integration and the increase in the dimensioning of the electrolyzers, it becomes crucial to model the electric dynamic behavior of the electrolyzer~\cite{EspinosaLpez2018}. We focus on polymer exchange membrane (PEM) electrolyzers and follow the modeling in~\cite{Lichtenberg22} to obtain a model that provides a dynamic description of the voltage and current dynamics as well as the outgoing hydrogen volumetric flow rates.

The electrical dynamics of the electrolyzer represented by the node $\n_{P,i}\in\mathcal N_P$ is essentially described by the activation overpotential $\mathbf{v}_{\text{a},i}$, which is a key component of the electrolyzer voltage that is mainly influenced by the reaction kinetics of the electrochemical process. It reflects the additional energy required to overcome the activation energy barrier for the hydrogen and oxygen evolution reactions.

The activation overpotential $\mathbf{v}_{\text{a},i}$ can be modeled by \cite{Lichtenberg22}  
\begin{equation}
C_{\text{a},i}\dot{\mathbf{v}}_{\text{a},i}= \mathbf{i}_{i}- R_{\text{a},i}^{-1}\mathbf{v}_{\text{a},i},
\label{eq:v_act}
\end{equation}
where $\mathbf{i}_{i}\in\mathbb R$ is the total electrolyzer current, $R_{\text{a},i}>0$ is the activation resistance, $C_{\text{a},i} = \tfrac{C_{\text{DL,cell,i}} A_{i}}{n_{c,i}}>0$ is the total double-layer capacitance, based on the cell surface area $A_{i}$ and $C_{\text{DL,cell,i}}$ is the cell capacitance and $n_{c,i}$ is the number of cells.

The total electrolyzer voltage $\mathbf{v}_{i}$ at the electrical node $\n_{P,i}\in\mathcal N_P$ is given as \cite{Lichtenberg22} 
\begin{equation}
\label{eq:v_ely}
\mathbf{v}_{i} = \mathbf{v}_{\text{oc},i} + \mathbf{v}_{\text{a},i} + \mathbf{v}_{\text{oh},i},
\end{equation}
with $\mathbf{v}_{\text{oc},i}$ being the open-circuit voltage of the electrolyzer, as described by the Nernst equation, i.e.,
\begin{align*}
    \mathbf{v}_{\text{oc},i}&=n_{c,i}\left( \mathbf{v}_{\text{st}} - \beta_T (T_i - T_{\text{st}}) + \frac{RT_i}{2F} \ln\left( \frac{p_{\text{H}_2}}{\sqrt{p_{\text{O}_2}} p_{\text{H}_2\text{O}}} \right)\right), 
\end{align*}
\noindent where $\beta_T = 0.0009 \, \text{V/K}$ is the temperature coefficient accounting for the variation of the open-circuit voltage with temperature $T_i>0$ and $F \approx 9.6485\cdot 10^4\, \text{C/mol}$ is the Faraday constant. The cell voltage under standard conditions $\mathbf{v}_{\text{st}}$ is typically $1.23 \, \text{V}$. Standard conditions are defined as a temperature of $ T_{\text{st}} = 298.15 \, \text{K} $ and a pressure of $ P_{\text{st}} = 101325 \, \text{Pa} $, and $p_{\text{H}_2}, p_{\text{O}_2}, p_{\text{H}_2\text{O}}$ are the partial pressures of hydrogen, oxygen, and water \cite{ESPINOSALOPEZ2018160}. 

The Ohmic overpotential $\mathbf{v}_{\text{oh},i}$ accounts for resistive losses in the electrolyte and electrodes and is given by 
\begin{align}
\mathbf{v}_{\text{oh},i}=n_{c,i}\frac{\delta_{\text{m},i}}{\sigma_{\text{m},i}A_{i}}\mathbf{i}_i,
    \end{align}
    with membrane thickness and conductivity $\delta_{\text{m},i},\sigma_{\text{m},i}>0$. 

The fuel cell can be thought of as the counterpart to an electrolyzer, i.e.\ hydrogen is consumed to produce electrical power.

Thus, the dynamics of the activation overpotential $\mathbf{v}_{\text{a},i}$ at the fuel cell node $\n_{P,i}\in\mathcal N_P$
in terms of the fuel cell current $\mathbf{i}_{i}$ is identical to the electrolyzer behavior \eqref{eq:v_act} 
and the total fuel cell voltage $\mathbf{v}_{i}$ is given by \cite{FuelCell} 
\begin{equation}
\label{eq:v_fuel}
-\mathbf{v}_{i} =-\mathbf{v}_{\text{oc},i} + \mathbf{v}_{\text{a},i} + \mathbf{v}_{\text{oh},i},
\end{equation}
with its output power depending on the hydrogen input flow rate. The volumetric flow rate of hydrogen can be expressed in terms of the electrolyzer and fuel cell currents $\mathbf{i}_{i}$ as \cite{CARMO2019165} 
\begin{equation}
\label{eq:q_ely}
\mathbf{q}_{\s,i} = \frac{M_{\text{H}_2}}{z\rho F}\mathbf{i}_{i},
\end{equation}
where $M_{\text{H}_2} = 2.016 \, \text{g/mol}$ is the molar mass of hydrogen, $z = 2$ is the number of electrons transferred per hydrogen molecule.

\begin{assumption}
\label{assum:ely}
Consider the electrolyzer and fuel cell dynamics \eqref{eq:v_act}, \eqref{eq:v_ely}, and  \eqref{eq:v_fuel}.
(i) The open circuit voltage $ \mathbf{v}_{\text{oc},i}$ is constant for all $i=1,\ldots,E+L$;
  
(ii) The activation resistance $R_{\text{a},i}$ is constant for all $i=1,\ldots,E+L$.
\end{assumption}

Assumption~\ref{assum:ely} is valid for moderate temperature changes within the electrolyzer and the fuel cell stacks, i.e., if the electrolyzers and fuel cells are operated relatively close to some desired operating state for the hydrogen volumetric flow rates and the temperature dynamics. This is reasonable in many settings, since often additional temperature control is applied to achieve a specific temperature set-point~\cite{Lichtenberg22,Pfennig25}. For larger temperature changes, the irreversible pH framework~\cite{RamL22} could be employed to derive a model that captures the associated thermodynamic more accurately.

Furthermore, a suitable constant value for the activation resistance $R_{\text{a},i}$ can be identified in experimental setups from measurement data~\cite{Yodwong21}. 

With Assumption~\ref{assum:ely}, we can rewrite the electrolyzer dynamics \eqref{eq:v_act} and \eqref{eq:v_ely} and the fuel cell dynamics \eqref{eq:v_fuel}, \eqref{eq:q_ely} as the following pH system with feed-through 
\begin{align}
    &\dot x_{\s,i}=(J_{\s,i}-R_{\s,i})\nabla H_{\s,i}(x_{\s,i})+B_{\s,i}u_{\s,i},\label{eq:ely_pH}\\ 
    &y_{\s,i}=B_{\s,i}^\top\nabla H_{\s,i}(x_{\s,i})+D_{\s,i}u_{\s,i}+d_{\s,i},
    \nonumber\\
    \nonumber
    &H_{\s,i}(x_{\s,i})=\tfrac12 C_{\text{a},i}^{-1}  x_{\s,i}^2, \quad x_{\s,i}=C_{\text{a},i}\mathbf{v}_{\text{a},i},\\
    &u_{\s,i}=\mathbf{q}_{\s,i},  \quad d_{\s,i}=\begin{cases}\frac{z\rho F}{M_{\text{H}_2}}\mathbf{v}_{\text{oc},i} &\!\!\!\! \text{if $1\leq i\leq E$,} \\ -\frac{z\rho F}{M_{\text{H}_2}}\mathbf{v}_{\text{oc},i} &\!\!\!\! \text{if $E+1\leq i\leq E+L$,}\!\! \end{cases}\nonumber\\
    & R_{\s,i}=R_{\text{a},i}^{-1}, \quad J_{\s,i}=0,\nonumber\\
    &B_{\s,i}= \frac{z \rho F}{M_{\text{H}_2}},\quad 
    D_{\s,i}=n_{c,i}\frac{\delta_{\text{m},i}}{\sigma_{\text{m},i}A_i}\frac{z^2\rho^2F^2}{M_{\text{H}_2}^2}
    \nonumber .
\end{align}
Note that the product of the electrolyzer input and output is equal to the electrical power 
\[y_{\s,i}u_{\s,i}=\frac{z\rho F}{M_{\text{H}_2}}\mathbf{v}_i\mathbf{q}_{\s,i}=\mathbf{i}_{i}\mathbf{v}_i,\] 
for $i=1,\ldots,E$, whereas the output of the fuel cell model is $y_{\s,i}=-\frac{z\rho F}{M_{\text{H}_2}}\mathbf{v}_{i}$ for all $i=E+1,\ldots,E+L$ reflecting the inverse relationship.

\section{Port-Hamiltonian modeling of integrated hydrogen systems}
\label{sec:network}

In this section, we interconnect the pH component models from Section~\ref{sec:pH_components} to obtain a pH model of the overall hydrogen grid. An example of such a grid is shown in Fig.~\ref{fig:HydrogenNetwork}.

The interconnected hydrogen grid without the sector coupling components is given by a combination of the pH models for the volumetric  flow rate dynamics for $\mathbf{q}$ that are given by \eqref{eq:pipe_pH}, \eqref{eq:comp_throttle_pH}, \eqref{eq:comp_flow_pH} and the pH models for the pressure dynamics for $\mathbf{p}$ that is given by \eqref{eq:storage_pH}, \eqref{eq:junction_pH}, \eqref{eq:comp_plenum_pH}. 

Thus, we define the state vector assigned to nodes and representing (weighted) pressures at storage units, junctions and compressors as
$$
x_\n=((x_{\n,i})_{i=1}^S, (x_{\n,i})_{i=S+1}^{S+F},(x_{\n,i})_{i=S+F+1}^{N}).
$$
Likewise, the state vector collecting the volumetric flow rates at the different pipes and compressors is defined as 
$$
x_\e=((x_{\e,i})_{i=1}^{M-2C},(x_{\e,i})_{i=M-2C+1}^M),
$$
where without loss of generality we group the edges $\e_f\in\mathcal{E}$ and $\e_m\in\mathcal{E}$ corresponding to any compressor $\n_i\in\mathcal N_C$, such that $m=f+1.$ Then, the overall system state vector is given by 
$
x_\g=(x_{\n},x_{\e}).
$

The pressure dynamics of all components of the hydrogen grid can be written as a diagonal combination of all subsystems~\eqref{eq:storage_pH}, \eqref{eq:junction_pH} and \eqref{eq:comp_plenum_pH}, i.e., 
\begin{align}
\label{eq:node_pH}
\dot x_\n&=(J_\n-R_\n)\nabla H_\n(x_\n)+B_\n u_\n,\\ 
y_\n&=B_\n^\top\nabla H_\n(x_\n)=\begin{bmatrix}y_\n^1\\ y_\n^2 \end{bmatrix}=\begin{bmatrix}\mathbf{p}\\ (\mathbf{p}_i)_{i=1}^{S+F} \end{bmatrix}, \nonumber \\
J_\n&=\diag(J_{\n,i})_{i=1}^N=0_{N\times N},\quad R_\n=\diag(R_{\n,i})_{i=1}^N, \nonumber \\ B_\n&=\begin{bmatrix} B_{\n}^1 & B_{\n}^2 \end{bmatrix}=
\begin{bmatrix} I_N & \begin{bmatrix}I_{S+F}\\ 0_{C\times(S+F)} \end{bmatrix} \end{bmatrix}, \nonumber  \\u_\n&=(u_\n^1,\mathbf{q}_\ex) :=(
    (u_{\n,i}^1)_{i=1}^N, (\mathbf{q}_{\ex,i})_{i=1}^{S+F}
).\!\! \nonumber 
\end{align}
Furthermore, we rearrange the entries of the input vector such that $u_\n^1$ collects all first entries of $u_{\n,i}$ for all $i=1,\ldots,N$, 
which are precisely the volumetric flow differences and the Hamiltonian is given by $H_\n=\sum_{i=1}^NH_{\n,i}$.

Likewise, the volumetric flow rate dynamics of all components of the hydrogen grid can be written as a diagonal combination of all subsystems~\eqref{eq:pipe_pH}, \eqref{eq:comp_flow_pH} and \eqref{eq:comp_throttle_pH}, i.e.,
\begin{align}
\label{eq:edge_pH}
\dot x_\e&=(J_\e-R_\e(x_\e))\nabla H_\e(x_\e)+B_\e u_\e+d_\e,\\ 
y_\e&=B_\e^\top\nabla H_\e(x_\e)=\begin{bmatrix}
    y_\e^1\\ y_\e^2
\end{bmatrix}=\begin{bmatrix}\mathbf{q} \\ (\mathbf{q}_{M+2(i-C)-1})_{i=1}^C\end{bmatrix},\nonumber \\
\!\!\!\!J_\e&\!=\!\diag(J_{\e,i})_{i=1}^M\!=\!0_{M\times M},\quad \!\!\!\!\!\!R_\e(x_\e)\!=\!\diag(R_{\e,i}(x_{\e,i}))_{i=1}^M, \nonumber \\ B_\e&=\begin{bmatrix} B_{\e}^1 & B_{\e}^2 \end{bmatrix}=\begin{bmatrix}I_M & \begin{bmatrix}
0_{(M-2C)\times C}\\\diag(\begin{bmatrix}
    1& 0
\end{bmatrix}^\top)_{i=1}^C
\end{bmatrix} \end{bmatrix},
\nonumber \\u_\e&=( u_\e^1,\Delta\mathbf{p}) :=(
    (u_{\e,i}^1)_{i=1}^M, (\Delta\mathbf{p}_i)_{i=1}^C), \nonumber 
\end{align}
and the Hamiltonian is given by $H_\e=\sum_{i=1}^MH_{\e,i}$.

The interconnection of the node dynamics~\eqref{eq:node_pH} and edge  dynamics~\eqref{eq:edge_pH} is established via
\[
u_\n^1=B_{\mathcal{G}}\mathbf{q}=B_{\mathcal{G}}y_\e^1,\quad u_\e^1=-B_{\mathcal{G}}^\top \mathbf{p}=-B_{\mathcal{G}}^\top y_\n^1,
\]
which leads to the following pH model for the interconnected hydrogen grid
\begin{align}
&    \dot{x}_{\g} = (J_{\g} - R_{\g}(x_{\g})) \nabla H_{\g}(x_{\g}) + B_{\g}u_{\g} + d_{\g}, \label{eq:grid_pH} \\
 &   y_{\g} = B_{\g}^\top \nabla H_{\g}(x_\g)=\begin{bmatrix}(\mathbf{p}_i)_{i=1}^{S+F} \nonumber \\ (\mathbf{q}_{M+2(i-C)-1})_{i=1}^C\end{bmatrix},\, u_\g=\begin{bmatrix}
       \mathbf{q}_{\ex} \\ \Delta\mathbf{p}
   \end{bmatrix},\\ 
    &H_{\g}(x_\g) = \sum_{i=1}^{N}H_{\n,i}(x_\n)+\sum_{i=1}^{M}H_{\e,i}(x_\e),\,  \nabla H_{\g}(x_{\g}) \!=\! \begin{bmatrix} \mathbf{p} \\ \mathbf{q} \end{bmatrix}, \nonumber \\
    &R_{\g}(x_\g) =\begin{bmatrix}
        R_\n &0_{N\times M}\\0_{M\times N}& R_\e(x_\e)
\end{bmatrix},\quad 
   J_{\g} = \begin{bmatrix} 0_{N\times N} & B_{\mathcal{G}} \\ -B_{\mathcal{G}}^\top & 0_{M\times M} \end{bmatrix},\nonumber \\    &B_\g=\begin{bmatrix}
       B_\n^2 & 0_{N\times C} \\ 0_{M\times(S+F)}& B_\e^2
   \end{bmatrix}, \quad d_{\g} = \begin{bmatrix} 0_{N}\\ \left(d_{\e,j}\right)_{j=1}^{M}\end{bmatrix} \nonumber.
\end{align}

In the following, we add $E\geq0$ electrolyzers which are connected to the storage units at $\n_1,\ldots,\n_{E}\in{\mathcal{N}_S}$ and $L\geq 0$ fuel cells connected to the storage units at $\n_{E{+1}},\ldots,\n_{E+L}\in{\mathcal{N}_S}$ by setting $\mathbf{q}_{\s}=(\mathbf{q}_{\s,i})_{i=1}^{E+L}=(\mathbf{q}_{\ex,i})_{i=1}^{E+L}$. 
For the resulting system, we keep the hydrogen volumetric flow rates $\mathbf{q}_{\s}$ as input variables. This leads to the following sector-coupled hydrogen system with state vector $x_{\sg}=(x_{\g},(x_{\s,i})_{i=1}^{E+L})$ 
\begin{align}
\nonumber
&\dot x_{\sg}=(J_{\sg}-R_{\sg}(x_{\sg}))\nabla H_{\sg}(x_{\sg})+B_{\sg}u_{\sg}+d_{\sg,1},\\
&y_{\sg}=B_{\sg}^\top\nabla H_{\sg}(x_{\sg}) + D_{\sg}u_{\sg} + d_{\sg,2}, \label{eq:gas_comp_elec_pH} \\ \nonumber 
&R_{\sg}(x_{\sg})=\begin{bmatrix}
    R_{\g}(x_{\g}) &0_{(N+M)\times(E+L)} \\ 0_{(E+L)\times (N+M)} & R_\s
\end{bmatrix},\\ &J_{\sg}=\begin{bmatrix}
    J_\g&0_{(N+M)\times(E+L)}\\0_{(E+L)\times (N+M)}&0_{(E+L)\times(E+L)}
\end{bmatrix}\nonumber, \\
&u_{\sg}=\begin{bmatrix}
    \mathbf{q}_{\s}\\ (\mathbf{q}_{\ex,i})_{i=E+L+1}^{S+F}\\ (\Delta\mathbf{p}_{i})_{i=1}^C
\end{bmatrix}, \nonumber \\  &B_{\sg}=
\begin{bmatrix}B_\g\\
   \begin{bmatrix}
         \tfrac{z \rho F}{M_{\text{H}_2}} I_{E+L}
         &0_{(E+L)\times(N-E-L)}
    \end{bmatrix}  
\end{bmatrix} , \nonumber\\ &d_{\sg,1}=(d_{\g},
   0_{E+L}   ),\quad d_{\sg,2}=(    (d_{\s,i})_{i=1}^{E+L}, 0_{N-E-L} 
    ),\!\! \nonumber\\   &D_{\sg}=\begin{bmatrix}\diag(D_{\s,i})_{i=1}^{E+L} &0_{(E+L)\times(N-E-L)}\\ 0_{(N-E-L)\times(E+L)}&0_{(N-E-L)\times(N-E-L)}\end{bmatrix}, 
\nonumber
\end{align}
and with the overall Hamiltonian $$H_{\sg}(x_\sg)=H_\g(x_\g)+\sum_{i=1}^{E+L}H_{\s,i}(x_{\s,i}).$$

\section{Conclusion}
In this work, we have presented a unified systematic port-Hamiltonian (pH) framework for modeling hydrogen systems, including storage units, electrolyzers, fuel cells, compressors, and pipelines. The inherent passivity properties of pH systems open the door to a structured analysis and controller synthesis for this emerging class of new energy systems. In future work, we will focus on further refining the presented modeling approach, especially by relaxing some of the assumptions made  and complementing the hydrogen system with an explicit pH representation of the electrical power system.

\bibliographystyle{ieeetr}
\bibliography{references}

\end{document}